# Extended five-term nonlinear drag model for a wide range of cylinder wakes

Osama A. Marzouk[1]

## Abstract

The unsteady variations in the near wake of a moving cylinder induce lift and drag forces on it, which are customarily normalized and expressed in terms of nondimensional lift and drag coefficients. While there are already several wake oscillator models for either a fixed or moving cylinder, special attention was given to modeling the lift coefficient for the case of a fixed cylinder or the case of a cylinder with one-degree-of-freedom motion in the cross-stream direction. When the drag coefficient is molded for a fixed or two-degree-of-freedom moving cylinder, a two-to-one frequency relationship (or quadratic coupling) between the drag and lift coefficients was assumed in the literature. However, we report situations of the excited wake of a vibrating cylinder, where such a modeling assumption fails to reproduce the actual pattern of the drag coefficient. We excite the wake of the cylinder by vibrating it harmonically in straight lines, and we then investigate the effect of this mechanical harmonic excitation on the lift-drag coupling using three tools for nonlinear dynamics analysis, namely, (1) time domain, (2) projection of the limit cycle, and (3) power spectra. We perform this analysis under different motion cases with manifested lift-drag coupling types that call on an extended universal drag model that accommodates such cases. Based on this, we propose a new reduced-order drag model with both linear and quadratic coupling terms to the lift as well as a mean component (thus, the proposed model consists of five terms). We verified the accuracy of the proposed reduced-order drag model by testing its ability to reproduce the time-dependent drag coefficient signals at a low Reynolds number of 300. These drag coefficient signals were obtained by applying direct numerical simulation (DNS) to integrate the two-dimensional incompressible Navier–Stokes equations governing the near wake fluid flow, with the aid of the finite difference method (FDM). The proposed drag model helps in extending the wake oscillators to more general cases of fluid–structure interaction (FSI) or vortex-induced vibration (VIV).

✉ Osama A. Marzouk
osama.m@uob.edu.om

[1] College of Engineering, University of Buraimi, Al Buraimi 512, Sultanate of Oman

**Keywords** Cylinder · Drag model · FSI · Lift coefficient · Nonlinear · Two-to-one · VIV · Wake oscillator

## Nomenclature (Greek letters first)

- $\eta$    Cross-stream displacement of the cylinder (expressed as a fraction of the cylinder diameter $D$)
- $\nu$    Kinematic viscosity of the fluid, $\nu = \mu / \rho$
- $\rho$    Density of the fluid
- $\mu$    Dynamic viscosity of the fluid
- $\Omega$    Nondimensional angular frequency of the applied exciting harmonic vibration of the cylinder
- $\omega_S$    Nondimensional Strouhal angular frequency (for a fixed cylinder), $\omega_S = 2\pi f_S$
- $\psi$    Phase angle (angle of a frequency component relative to $a_{1L}$)
- $a_{1L}$    Frequency component of $C_L$ at $\omega_S$ (or $\Omega$)
- $a_{2L}$    Frequency component of $C_L$ at $2\omega_S$ (or $2\Omega$)
- $a_{3L}$    Frequency component of $C_L$ at $3\omega_S$ (or $3\Omega$)
- $a_{1D}$    Frequency component of $C_D$ at $\omega_S$ (or $\Omega$)
- $a_{2D}$    Frequency component of $C_D$ at $2\omega_S$ (or $2\Omega$)
- $C_D$    Drag coefficient (time-dependent)
- $C_f$    Friction coefficient (time-dependent)
- $C_L$    Lift coefficient (time-dependent)
- $C_P$    Pressure coefficient (time-dependent)
- $D$    Diameter of the cylinder
- $L$    Length of the cylinder (L = 1 for the present two-dimensional flows)
- $f_S$    Nondimensional Strouhal cyclic frequency (for a fixed cylinder)
- $p$    Pressure
- Re    Reynolds number, Re $= \rho U D / \mu = U D / \nu$
- $U$    Free-stream velocity

## 1 Introduction

The flow over a cylindrical body is commonly encountered in different applications [1–3], such as offshore oil risers [4, 5], wind turbine towers [6–8], solar updraft systems [9, 10], overhead power distribution poles [11–13], chimneys of power plant and combustion processes [14–16], above-ground piping network of gases or liquids [17–19], and various structural elements in the buildings sector or the energy sector [20–23]. Therefore, theoretical investigations of such a problem have been conducted by several researchers [24–26].

There are several reduced-order models (ROMs) [27, 28] for the lift and drag coefficients for a fixed cylinder. These lift and drag coefficients are nondimensionalized forms of the two orthogonal components of the exerted hydrodynamic force due to the time-dependent variations in the near wake [29]. As in other fluid dynamics or gas dynamics applications [30, 31], the lift is the force component acting on the cylinder in the cross-stream direction, while the drag is the force component acting on the cylinder in the free-stream direction [32, 33]. For convenience and simplified discussion in the current work, we may use the term "lift" or "lift force" to refer to the lift coefficient [34], and use the term "drag" or "drag force" to refer to the drag coefficient [35]. This should not pose a fundamental concern since either force coefficient is derived from its relevant force component using the same scaling factor [36–38], which is $\frac{1}{2} \rho U^2 D L$, where $\frac{1}{2} \rho U^2$ is the dynamic pressure, and $D L$ is the projected area of the cylinder perpendicular to the incoming flow [39, 40]. While some of the recent phenomenological reduced-order models (i.e., models that explain one or more observed phenomena) and fully analytical reduced-order models (i.e., models whose parameters are fully related to some variables of the flow case of interest and thus can be systematically estimated) have proved to be very accurate and can even model the transient phase, in addition to the steady-state phase of development of these forces very well [41], they need to be extended to the moving-cylinder case. In fact, some early attempts to do that go back to 1970 [42], and there are existing models that represent the induced lift force for a moving cylinder in a free stream direction (in-line motion) [43, 44]. Several of these models focus on the case of a cylinder that moves in the cross-stream direction [45, 46], and these models couple this motion to the lift oscillator, which is a nonlinear ordinary differential equation (ODE) representing a self-excited dynamical system [47–49], without paying attention to the drag [50–58].

Other models that consider the drag force as well as the lift impose a specific coupling between this drag force and the lift force, which is purely quadratic [59]. Another example is the work based on modeling the structural response of elastic cable suspensions, where a model of the drag was introduced, and the drag was represented as a function of time and the cable's longitudinal coordinate in the static condition [60]. In that work, the drag was decomposed into the product of a temporal (time-dependent) function and a spatial (space-dependent) function. The time-dependent function is modeled by a van der Pol equation (cubic nonlinear oscillator) whose linear natural frequency is twice the shedding frequency (which was assumed to be equal to the angular natural frequency of the cable). As in drag-modeling studies for a fixed cylinder, that research work used quadratic coupling between the lift (which is modeled by another van der Pol equation whose linear natural frequency is equal to the shedding frequency) and the drag. Therefore, the two-to-one frequency relationship between the lift and the drag (the fundamental frequency of the oscillatory drag signal is twice the fundamental frequency of the oscillatory lift signal) was assumed in that work. While this assumption is usually acceptable for a wide range of recorded data of VIV (vortex-induced vibration), we show here in the current study that this assumption is not accurate for some types of flows, such as when the cylinder is vibrating at a prescribed tilt angle relative to the free stream [61, 62]. We denote such cases when the two-to-one frequency relationship between the lift and drag is not satisfied as "non-traditional wakes"; this is in contrast to the case when this relationship is satisfied, which we denote by "traditional wakes". Also, the discussed semi-empirical drag model [60] included two empirical constants that were selected and were not well related to the data of interest being modeled (the training dataset). The proposed reduced-order drag model here is fully analytical, meaning that we provide closed-form mathematical expressions to relate its parameters to some variables that are extracted from the training data to be modeled such that good agreement between the reduced-order model (ROM) results and these training data is achieved.

A research team also worked on modeling both the lift and drag forces for fixed and moving cylinders with two degrees of freedom [63]. Instead of using a wake oscillator, they followed a force-decomposition approach, where the lift was assumed to be composed of a main frequency component at the shedding frequency and the odd harmonics, whereas the drag was assumed to be composed of a mean value in addition to a main frequency component at twice the shedding frequency and the even harmonics. For a fixed cylinder, the harmonics were neglected, and one-term lift and one-term drag (in addition to the mean-drag value) were used to represent the signals (time series) of these two fluid forces, which are used to excite the cylinder and move it in both the $x$ and $y$ directions. An earlier reduced-order model (ROM) [41] accounted for the third harmonic of the lift when modeling the lift induced on a fixed cylinder. That drag model is consistent with the one used by another research team [63] in terms of using one harmonic only for the drag of a fixed cylinder. For a moving cylinder, the work described in [63] claimed that the harmonics can still be dropped if the shedding frequency is away from the natural frequency of the cylinder, but these harmonics should be retained if the two frequencies are close to each other (i.e., at resonance [64, 65]). The coefficients of these harmonics can be obtained from any arbitrary available time-history lift and drag data using auto-regressive moving averaging (ARMA) [66–68]. In the work described in [63], it was still assumed that a two-to-one frequency relationship between the lift and drag exists, and even harmonics in the lift as well as odd harmonics in the drag were totally ignored. Again, these assumptions lead to a model that fails to model the drag associated with a cylinder undergoing a harmonic inclined vibration, as we demonstrate here.

We contribute to this class of problems by extending existing drag reduced-order models (ROMs) that are based on a simple quadratic relationship between the lift and drag, to situations where this relationship does not hold. One of these situations is demonstrated here, which occurs when the cylinder vibrates along directions other than the cross-stream. This can be observed in fluid–structure interaction (FSI) or vortex-induced vibration (VIV) problems, where the oscillatory motion of the cylinder (or cylinder-like element) cannot be restrained to the perpendicular cross-stream direction [69–71].

Figure 1 illustrates this generalized situation, where the cylinder undergoes an oscillation that is inclined (i.e., lies between the in-line direction and the cross-stream direction).

The change in the lift-drag relationship is accompanied by other changes in their spectral characteristics. However, these changes are significant for the drag, and thus the available models fail in these cases. When the traditional wake regime and its implied simple two-to-one lift-drag relationship and spectra are restored, the proposed reduced-order model (ROM) automatically converges to a form that is identical to traditional existing drag models (the extra terms we add in the proposed drag ROM automatically vanish).

**Fig. 1** Illustration of the cylinder motion and surrounding fluid

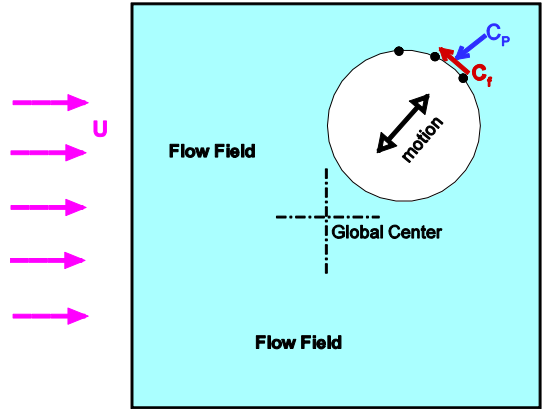

## 2 Flow analysis

In order to provide a training dataset of fluid force signals with a drag-lift relationship that does not obey the traditional two-to-one frequency relationship (we need such a dataset to develop and test the new universal reduced-order drag model), and compare its performance to existing drag models that are limited to the case of two-to-one frequency relationship between the lift and drag; we generate this training dataset using computational fluid dynamics (CFD) modeling [72–75]. For this part, we numerically solve (integrate in time and space) the two-dimensional incompressible Navier–Stokes (NS) Eqs. [76–78], which are

$$\frac{\partial u}{\partial x} + \frac{\partial v}{\partial y} = 0$$
$$\frac{\partial u}{\partial t} + \frac{\partial u^2}{\partial x} + \frac{\partial uv}{\partial y} = -\frac{1}{\rho}\frac{\partial p}{\partial x} + \nu\frac{\partial^2 u}{\partial x^2} + \nu\frac{\partial^2 u}{\partial y^2} \quad (1)$$
$$\frac{\partial v}{\partial t} + \frac{\partial uv}{\partial x} + \frac{\partial v^2}{\partial y} = -\frac{1}{\rho}\frac{\partial p}{\partial y} + \nu\frac{\partial^2 v}{\partial x^2} + \nu\frac{\partial^2 v}{\partial y^2}$$

In the above equation, $u$ and $v$ are the axial (streamwise) and vertical (cross-flow) components of the fluid flow velocity vector in the $x$ and $y$ coordinates, respectively. The symbol $t$ designates the time, $\rho$ is the fluid density, $p$ is the fluid pressure, and $\nu$ is the fluid's kinematic viscosity. Due to the lack of heat sources or sinks [79], heat transfer mechanisms, thermal radiation [80], temperature-dependent fluid properties [81], and high-speed compressible-fluid viscous heating [82], the energy equation is not solved [83, 84]. This implies a reasonable assumption of isothermal (constant-temperature) cold flow.

Before the numerical solution, these NS equations are nondimensionalized using the free-stream velocity $U$, the cylinder diameter $D$, and twice the dynamic pressure $\rho U^2$. We then discretize and solve the resulting system of partial differential equations (PDEs) in the time domain using direct numerical simulation (DNS) [85, 86]. We then compute the unsteady flow field around the vibrating cylinder, including the pressure distribution and the shear-stress distribution on its surface.

**Fig. 2** Illustration of the pressure and friction coefficients at the cylinder surface

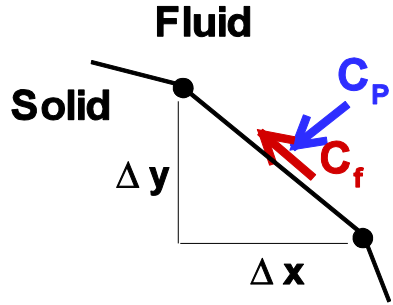

Although the lift coefficient ($C_L$) and the drag coefficient ($C_D$) for a cylinder should be theoretically obtained through integrating the nondimensional forces over the cylinder's surface (over circular circumference in the current two-dimensional problem), in computational fluid dynamics (CFD) modeling [87], such continuous symbolic integration needs to be discretized [88, 89] and replaced by an quantitative summation process [90]. This transformation is needed because a closed-form analytical solution does not exist for the lift coefficient or the drag coefficient for a cylinder subject to a viscous fluid flow [91, 92], although it exists in a very simplified case of potential flows where the viscosity is totally neglected [93, 94]. Thus, numerical approximation is necessary, and it is achieved by representing the perfectly circular circumference of the cylinder as a series of connected straight-line segments. Each segment is formed by connecting two adjacent grid points located on the cylinder's surface (circumference). The incremental distances $\Delta x$ and $\Delta y$ are the horizontal and vertical projections of the distance between two adjacent grid points on the cylinder surface. Figure 2 illustrates a portion of the cylinder surface between two grid points and shows $\Delta x$ and $\Delta y$.

We numerically integrate the nondimensional pressure coefficient and the nondimensional friction coefficient over the cylinder surface to compute the corresponding nondimensional lift and drag coefficients, which are recorded as two time signals according to

$$\begin{aligned} C_L(t) &= \oint C_P(t)dx + \oint C_f(t)dy \\ C_D(t) &= -\oint C_P(t)dy + \oint C_f(t)dx \end{aligned} \quad (2)$$

where $C_P$ is the time-dependent pressure coefficient (or nondimensional pressure). It is the gauge pressure (the absolute pressure minus the freestream ambient or atmospheric pressure) divided by the dynamic pressure as a scaling factor [95–97]. Similarly, $C_f$ is the time-dependent friction coefficient (or nondimensional shear stress). It is the shear stress divided by the dynamic pressure as a scaling factor.

At the present Reynolds number (Re) of 300, the pressure contribution (the first term on the right-hand side of Eq. 2) to the lift and drag coefficients is more significant than the friction contribution (the second term on the right-hand side of Eq. 2). However, we did not neglect the friction term when computing the lift and drag coefficients.

It is useful to mention here that the pressure contribution to the lift and drag coefficients arises from the asymmetric pressure distribution over the cylinder's surface [98]. Vertical asymmetry (deviation between the pressure distribution over the top half of the cylinder and the bottom half of the cylinder) contributes to the lift coefficient, while horizontal asymmetry (deviation between the pressure distribution over the front, upstream half of the cylinder and the rear, downstream half of the cylinder) contributes to the drag coefficient. This pressure asymmetry, in turn, is primarily caused by the low-pressure wake zone formed downstream of the cylinder due to eddies (vortices) [99]. A vortex has a low-pressure zone at its core [100] due to the high spinning speeds developed there, in accordance with Bernoulli's principle (energy conservation) [101], where pressure energy is converted into kinetic energy, as well as due to the centrifugal effects that pull the spinning fluid outward and establish a suction power. It is also useful to elaborate on the flow regimes over a stationary cylinder based on the Reynolds number. At extremely low Reynolds numbers ($Re < 5$), reflected by a very slow incoming flow, the flow can be described as "creeping" [102]. For such a limiting case, the flow around the cylinder is nearly double-symmetric (having both front-rear symmetry as well as top–bottom symmetry), and no separation of the boundary layer occurs. The creeping flow is laminar and steady, where no irregularities or unsteadiness are observed. This nearly double-symmetric flow regime and double-symmetric pressure distribution cause a cancellation of the pressure effects, leading to vanishing contributions to either the lift coefficient or the drag coefficient. However, the skin friction effect still results in a non-zero drag coefficient. As the Reynolds number increases up to around 40, the boundary layer separates at the downstream side of the cylinder's surface, and two small vortices (Fopple vortices) develop [103, 104], but they remain confined within the near-wake zone rather than being shed away downstream. Thus, this flow regime features a static (steady) separated near-wake region. As the Reynolds number increases further up to around 200, a laminar periodic vortex street is established in the wake. This oscillatory pattern of the flow around the cylinder in its near wake zone causes the lift and drag coefficients to also exhibit an oscillatory behavior. The laminar vortex street gradually transitions into a turbulent one until a fully turbulent vortex street is reached, approximately at a transition threshold Reynolds number or critical Reynolds number [105] of $3 \times 10^5$ [106]. Up to this Reynolds number, the boundary layer over the unseparated part of the cylinder's surface is laminar. Laminar boundary layers exhibit relatively weak resistance to separation compared to turbulent boundary layers. This causes the laminar boundary layer to separate early [107], and this enlarges the width of the low-pressure wake. This magnifies the imbalance in the oscillatory pressure distribution, leading to a large contribution to the drag and the lift coefficients. The portion of the drag caused by the pressure imbalance is called form drag [108] or pressure drag [109]. On the other hand, the skin friction effect of a laminar boundary layer is smaller than the skin friction effect of a turbulent boundary layer. This makes the frictional contribution to the drag small compared to the pressure when the Reynolds number is below the transition threshold. As the Reynolds number increases beyond this transition threshold, the boundary layer undergoes a gradual change from a laminar type to a turbulent type, and eventually becomes a turbulent boundary layer. This retards the boundary layer separation, causing the two separation points to move further downstream along the cylinder's surface, leading to a narrower zone of low-pressure wake. This suppresses the pressure imbalance, and this weakens the influence of the pressure distribution on the lift and drag coefficients. On the other hand, the larger skin friction influence of the turbulent boundary layer increases its contribution to the lift and drag coefficients.

The cylinder undergoes a sinusoidal vibration with a frequency that is equal to the Strouhal frequency (i.e., $\Omega = \omega_S$), and with a prescribed vibrational amplitude. The cylinder's vibrational motion starts from the beginning of the CFD simulation.

The values of $C_P$ and $C_f$ are resolved at the grid points; then, an average value is considered to be applied at the midpoint. These midpoint values are used in Eq. 2 to find $C_L$ and $C_D$ at every time step after the flow field is resolved. The arbitrary Lagrange Euler (ALE) technique was used to implement the cylinder motion and couple it to the flow solver [110–113]. In this ALE technique, Eq. 1 is written in an inertial coordinate system (having a global fixed center). Initially, the local cylinder center coincides with this global center. The grid deforms every time step, and it is re-generated algebraically to fit the new fluid region.

A rigorous numerical approach was followed to ensure that the resolved flow features are physical and do not reflect numerical errors. Several test cases were performed, which include solving the flow field past a fixed cylinder and a moving cylinder under different simulation parameters (e.g., time step, grid resolution, and grid size) so that we get converged solutions that are satisfactorily independent of these preset parameters that are manually selected rather than automatically derived.

In order to validate the adopted numerical simulation procedure, we performed several calculations in order to reproduce other reported benchmarking results [114–116] from numerical or experimental studies for both fixed and moving cylinders. For example, we ran a test case at Reynolds number 500 for both fixed and moving cylinders and compared the results with a similar two-dimensional direct numerical simulation of an earlier independent study [117], where they used the spectral element method with ninth-order Gauss–Lobatto-Legendre polynomials [118]. For the fixed cylinder, we got $f_S = 0.225$ (compared to 0.228), time-averaged $C_D = 1.429$ (compared to 1.46), and an amplitude of $C_L = 1.159$ (compared to 1.20). Our values for this case also agree with the measured value of $f_S = 0.22$ [119] and the measured value of time-averaged $C_D = 1.5$ [120]. We also solved the case of a moving cylinder in the cross-stream direction at different frequencies and got similar results as reported by our benchmarking study [117]. In Fig. 3, we show the agreement of our results with the reported ones in terms of the limit cycle projection onto the $C_L$–$\eta$ plane at motion frequencies $0.875 f_S$ and $0.975 f_S$. The lift response is periodic in both cases.

## 3 Model analysis

We start this section by illustrating the quadratic coupling between the drag and lift for the fixed cylinder. Figure 4 shows typical power spectra of the lift and drag coefficients for this case, which we obtained from the numerical simulation. The main component of the lift occurs at $f_S$. We denote this component by $a_{1L}$. The main component of the drag occurs at twice this frequency; that is why we denote it by $a_{2D}$ (the numeral "2" in the subscript means twice the frequency of the fundamental frequency component of the lift). As a consequence of this two-to-one frequency relationship, the projection of the limit cycle of the wake in the lift-drag plane consists of two loops, as shown in Fig. 5. This two-loop pattern was also observed at a higher Reynolds number of 500 [121].

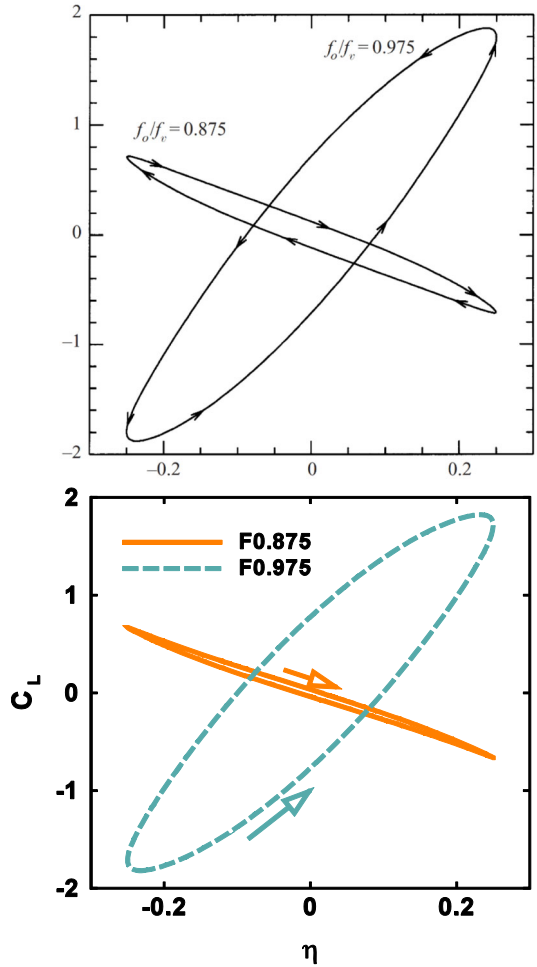

**Fig. 3** Projection of the limit cycle onto the $C_L$-$\eta$ plane for a moving cylinder in the cross-stream direction. Results from our simulations (bottom) are compared with those from the simulations reported in [117] (top), reproduced with permission (License number 6026501366615 dated 12/May/2025, obtained through the CCC "Copyright Clearance Center" RighsLink® service)

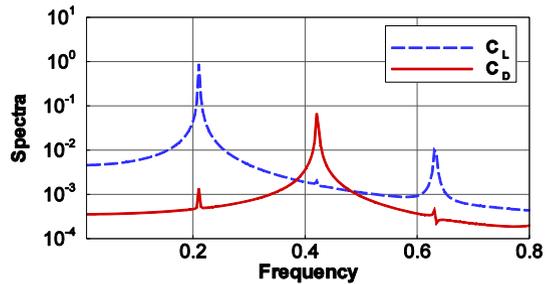

**Fig. 4** Power spectra of the lift and drag for the fixed cylinder

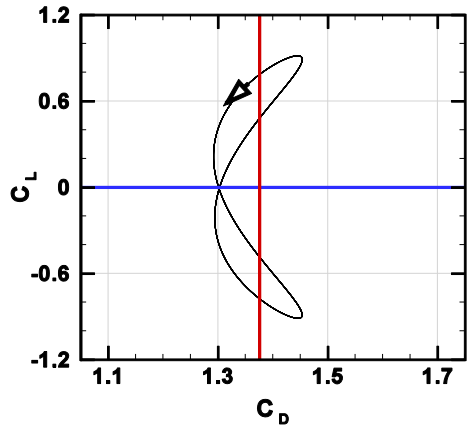

**Fig. 5** Projection of the limit cycle onto the $C_L$-$C_D$ plane for the fixed cylinder

Based on this two-to-one frequency relationship, proposed drag models for the fixed cylinder use quadratic coupling terms with the lift [122, 123]. The first proposed model for this category was

$$C_D = C_{D,\text{ mean}} - 2\,\frac{a_{2D}}{\omega_S\, a_{1L}^2}\, C_L\, \dot{C}_L \quad (3)$$

The rationale for the above ROM drag model is that it has the least number of necessary terms, which is two. So, it was made as simple as possible as an initial attempt to model the lift-drag coupling. The above ROM is based on the idea of splitting the drag coefficient signal into two independent components. First, there is a constant term ($C_{D,\text{mean}}$), which is the mean (the time-averaged value or the DC component) of the drag coefficient signal. It is represented by the first term in the above ROM equation. Then, there is the oscillatory component (AC component) of the drag coefficient signal, which has a zero mean. It is represented by the second term in the above ROM equation. Given that the lift and drag are Cartesian force components of the same hydrodynamic force exerted on the cylinder surface, it was considered reasonable in the literature to relate the oscillatory component of the drag to the lift, which is a pure oscillatory signal (has zero mean). In addition, because a two-to-one frequency relationship was observed in the literature between the oscillatory component of the drag and the oscillatory lift, it was decided to model the oscillatory component of the drag as a quadratic function of the oscillatory lift signal. This quadratic dependence can take six forms, namely $(C_L^2), (-C_L^2), (\dot{C}_L^2), (-\dot{C}_L^2), (\dot{C}_L C_L)$, and $(-\dot{C}_L C_L)$. Because the phase angle between the oscillatory component of the drag and the oscillatory lift was found in the literature to be near 270°, the last option was considered the most suitable. The multiplier factor of $2\frac{a_{2D}}{\omega_S a_{1L}^2}$ in the second term of the above two-term

ROM drag model is added for proper scaling, such that the reconstructed time-domain drag coefficient signal from the ROM equation is close in magnitude to the observed training data for a fixed cylinder.

While the above drag model is able to reproduce the simulated $C_D$ in the time domain to a good extent, it cannot reproduce the phase angle between the drag and the lift (i.e., the phase with which $a_{2D}$ leads $a_{1L}$). Therefore, the ROM represented by Eq. 3 was extended to another amended ROM [41] that has two quadratic terms, as shown in Eq. 4, instead of one term. The weights $q_1$ and $q_2$ are related to the phase angle found in the training dataset. This amended ROM is

$$C_D = C_{D,\text{ mean}} + 2q_1 \frac{a_{2D}}{a_{1L}^2}\left(C_L^2 - \left[C_L^2\right]_{\text{mean}}\right) + 2q_2 \frac{a_{2D}}{\omega_S\, a_{1L}^2} C_L \dot{C}_L \quad (4)$$

In the current study, we aim to extend this amended ROM further such that it holds valid for fixed and moving cylinders with a more general lift-drag relationship than the commonly assumed quadratic one. In the latter case, the two-to-one frequency relationship is distorted. To identify the behavior of the appropriate coupling in these new situations, we analyzed the simulated lift and drag coefficients for a vibrating cylinder at a vibration amplitude of 7.5% $D$. This value is relatively small for typical problems of VIV (vortex-induced vibration), but we selected it after conducting an extensive amount of CFD simulations at lower and higher amplitudes. Our amplitude selection criterion is the preference for an amplitude value that manifests a wide variety of lift and drag patterns as the motion inclination angle of the cylinder is varied while keeping the amplitude of this motion unchanged. This is important to achieve the aims of the current study and to provide challenging profiles of these forces (i.e., the lift and drag) and their relationship. Lower amplitudes can result in force patterns with a two-to-one relationship (or close). Existing drag models for this type already perform satisfactorily, and we then become unable to demonstrate the usefulness of the proposed new model. On the other hand, higher oscillation amplitudes may result in the dominance of chaotic or quasi-periodic force patterns. In these cases, the limit cycle is highly distorted, and we then become unable to demonstrate the relationship between the model parameters and the temporal and spectral properties of the lift and drag. The reader is reminded that this study is not concerned primarily with the CFD simulations and numerically resolving the flow field or the lift and drag forces, although it is a large component of the current work. Rather, the CFD simulation conducted here shows how existing reduced-order drag models need to be extended to accommodate cases when the lift and drag relationship does not follow a commonly perceived type for which these models were derived. Also, the CFD simulation provides training datasets for computing the proposed reduced-order model (ROM) for the drag coefficient. We later show several analysis techniques that lead to the proper modifications of these existing reduced-order drag models, and emphasize the influence of the proposed modifications on the modeled results. The proposed reduced-order drag model here contains some parameters that can then be tuned (based on closed-form rigorous analytical expressions, not by trial and error) to other training datasets generated at higher motion amplitudes than the selected 7.5% $D$. The Reynolds number is maintained at 300, which is adequate here.

We first look at the projection of the limit cycle onto the $C_L$-$C_D$ plane and the distortion in the two symmetric loops about the $C_L = 0$ axis. In Fig. 6, the motion axis is rotating toward the cross-stream direction as we proceed from case (a) to case (f). The inclination angles of the oscillation line are from 40° in case (a), 50° in case (b), 60° in case (c), 70° in case (d), 80° in case (e), and 90° (thus cross-stream oscillation) in case (f). We show the change of the lift-drag relationship in the projection of the limit cycle onto the plane of their nondimensional coefficients ($C_L$ and $C_D$) through this figure. The two-to-one frequency relationship (or the two-loop curve) that existed for the fixed-cylinder case is altered for all cases except for the case of cross-stream motion (the traditional case of one-degree-of-freedom motion). As the motion axis tilts away from the cross-stream direction, the lower loop (negative lift) becomes thinner while the upper loop (positive lift) becomes thicker. Eventually, the lower loop disappears, and we end up with one loop only, as in case (a) within Fig. 6. This fact suggests adding a linear coupling of the drag to the lift. In case the two-loop curve is transformed into a multi-loop curve, then a better alternative is to propose a higher-order coupling (e.g., cubic) rather than a lower-order coupling (i.e., linear). This simple tool of graphical analysis for this nonlinear dynamical system revealed a lot of useful details about the dynamics of the problem, and allowed qualitative determination of the appropriate modeling type. However, this elementary analysis should be complemented with a more thorough and quantitative one. For this purpose, we analyze the same lift and drag signals (time series) in the spectral domain and identify the frequency components of each signal. The power spectra of the lift and drag for the cases presented in Fig. 6 are shown in Fig. 7. Again, case (f) corresponds to a cross-stream oscillatory motion. The spectra of the lift and drag with such a cross-stream motion are qualitatively similar to those obtained for the fixed cylinder (was shown in Fig. 4). This confirms what we mentioned earlier about the similarity of the two cases when examined in the lift-drag plane. So, we were able to expose new lift-drag relationship patterns that were not identified before. With a cross-stream motion, the lift is composed of a dominant fundamental lift frequency component ($a_{1L}$) and its odd harmonics (such as $a_{3L}$), whereas the drag is composed of a mean quantity ($C_{D, \text{mean}}$) plus a dominant drag frequency component ($a_{2D}$) and its even harmonics (such as $a_{4D}$). So, $a_{2L}$ and $a_{1D}$ in this case are very small and can be neglected. As the motion axis tilts away from the cross-stream direction, $a_{1L}$ increases slowly, and $a_{3L}$ and $a_{2D}$ show weak non-monotonic variations. However, $a_{2L}$ and $a_{1D}$ increase rapidly. Since $a_{2L}$ remains about one order of magnitude less than $a_{1L}$, the distortion of the periodic lift signal is small. In contrast, $a_{1D}$ approaches and even exceeds $a_{2D}$ at some point, which is reflected in the strong deviation in the drag signal from the periodic pattern obtained either with a cross-stream motion or without any motion (a fixed cylinder). Based on the aforementioned findings and analysis, we propose the following extended (we refer to it as "universal") reduced-order drag model with five terms:

$$C_D = C_{D,\text{mean}} + \lambda_1 \frac{a_{2D}}{a_{1L}} C_L + \lambda_2 \frac{a_{2D}}{\Omega\, a_{1L}} \dot{C}_L \\ + 2q_1 \frac{a_{2D}}{a_{1L}^2} \left( C_L^2 - \left[ C_L^2 \right]_{\text{mean}} \right) + 2q_2 \frac{a_{2D}}{\Omega\, a_{1L}^2} C_L \dot{C}_L \quad (5)$$

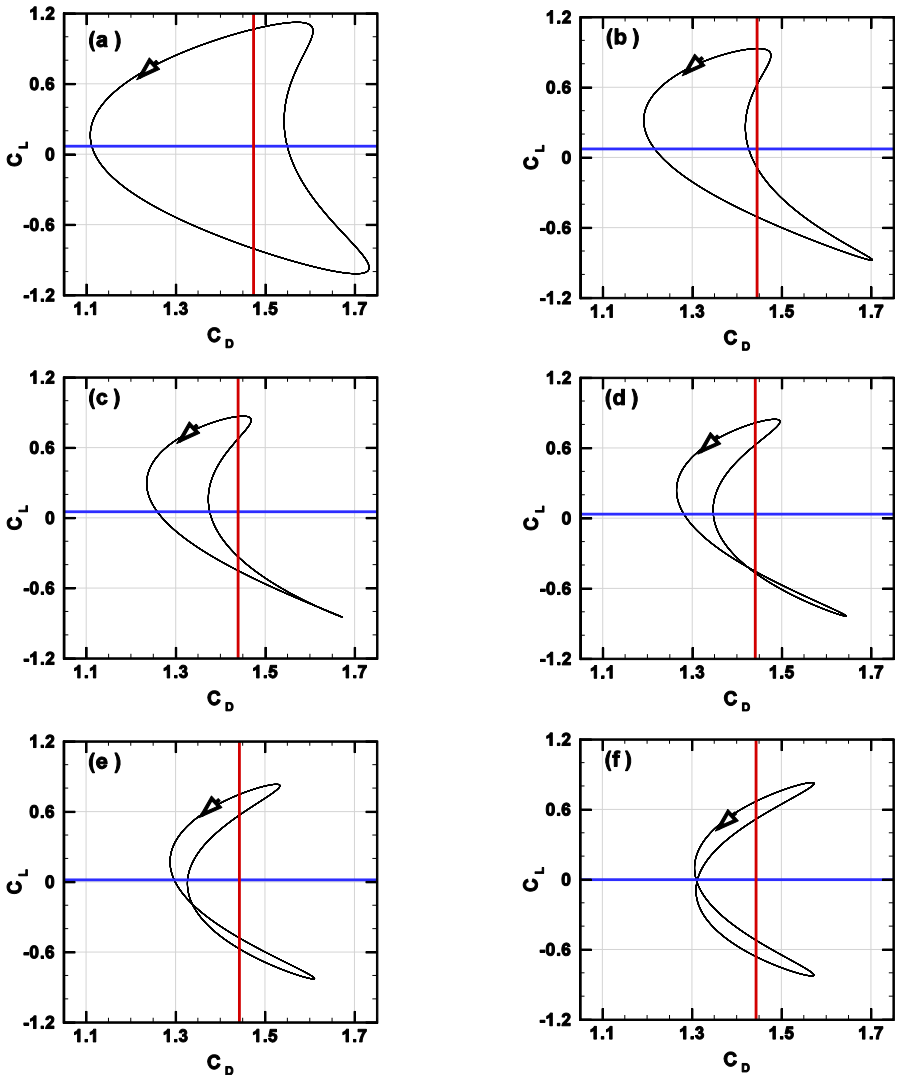

**Fig. 6** Projection of the limit cycle onto the $C_L$-$C_D$ plane for a moving cylinder compared at different cases with the motion axis is approaching the cross-stream direction. The motion is exactly along the cross-stream direction in case (f)

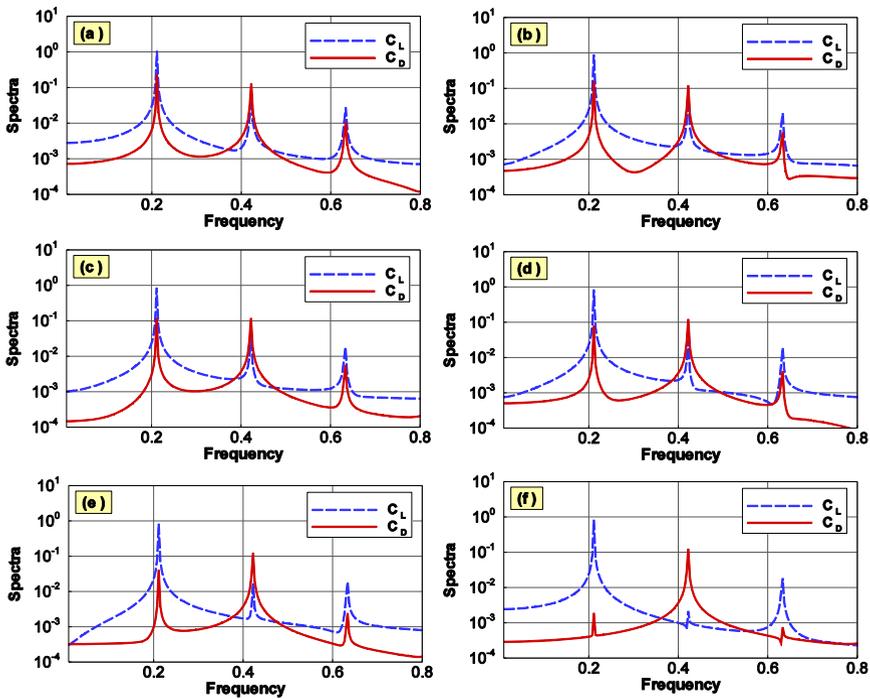

**Fig. 7** Spectra of the lift and drag coefficients for the moving cylinder cases in Fig. 6

The quadratic terms having the model coefficients $q_1$ and $q_2$ are almost identical to those presented earlier in Eq. 4 for the amended drag ROM, but the Strouhal frequency $\omega_S$ in Eq. 4 is replaced here by the sinusoidal vibration frequency $\Omega$ (both are equal in the current work). The new proposed linear-term coefficients $\lambda_1$ and $\lambda_2$ are determined such that the phase of $a_{1D}$ relative to $a_{1L}$ (we denote this phase angle by $\psi[a_{1D}, a_{1L}]$) is matched with the training CFD-based value. So

$$\lambda_1 = \cos(\psi[a_{1D}/a_{1L}]); \quad \lambda_2 = \sin(\psi[a_{1D}/a_{1L}]) \quad (6)$$

Similarly, the existing coefficients $q_1$ and $q_2$ for the quadratic terms are determined such that the phase of $a_{2D}$ relative to $a_{1L}$ (we denote this phase angle by $\psi[a_{2D}, a_{1L}]$) is matched with the training CFD-based value. So

$$q_1 = \cos(\psi[a_{2D}/a_{1L}]); q_2 = \sin(\psi[a_{2D}/a_{1L}]) \quad (7)$$

In the old reduced-order drag models of Eq. 3 (the primitive two-term drag ROM) and Eq. 4 (the amended three-term drag ROM), the angular shedding frequency $\omega_S$ appears because for a fixed cylinder [124] (which was the scope of these models), the frequency of the main lift component takes place at this value. However, in the proposed five-term ROM of Eq. 5, we replaced this angular shedding frequency $\omega_S$ with the angular vibration frequency $\Omega$, which can be different from $\omega_S$ (although $\Omega$ is set to be equal to $\omega_S$ here). This is to account for the fact that the shedding frequency (and the frequency of the main lift component) when the cylinder vibrates corresponds to the vibration frequency rather than to the Strouhal frequency that pertains to the fixed cylinder, when the wake is not excited by the mechanical vibration.

## 4 Results

We solve the flow field around a vibrating cylinder by direct numerical simulation (DNS) of the two-dimensional Navier–Stokes equations at Re 300, which is low enough to justify the two-dimensionality and the use of DNS without turbulence modeling [125–130]. It should be noted that the previous statement does not imply a connection between the dimensionality of the Navier–Stokes equations and the Reynolds number. For example, turbulence (at high Reynolds number) can be investigated using the two-dimensional Navier–Stokes Eqs. [131]. On the other hand, laminar flow dynamics (at low Reynolds numbers) can be studied using the three-dimensional Navier–Stokes Eqs. [132, 133]. Instead, we point out that in a real flow regime at Re 300, the three-dimensionality and turbulence exist, but not at a predominant scale [134, 135]. The transition from a two-dimensional to a three-dimensional flow field, and from a laminar to a turbulent flow field, happens gradually as the Reynolds number increases. The selected Reynolds number of 300 is viewed as a suitable trade-off between (1) much lower values in which the flow is trivially steady without interesting dynamics, and (2) much larger values in which the complexities of the real flow are partly missed due to the inability of the two-dimensional laminar version of the Navier–Stokes equations to capture them.

The time-dependent lift and drag coefficients are recorded as time signals during each case of computational fluid dynamics (CFD) simulation, and these time signals are obtained after a stable periodic time response is reached.

High-order spectral analysis is then applied to these signals to evaluate the significant components and their lead angles (relative to $a_{1L}$). These evaluated spectral parameters of the simulated lift coefficient are used in Eqs. 5–6 to generate a fitted drag-coefficient signal (i.e., obtained by applying the proposed five-term reduced-order drag model after identifying the numerical values of its five parameters through closed-form analytical expressions [136, 137]), which is then compared with the one obtained from the CFD simulation (the training set). We conducted several cases and got a good agreement in them. Some of these cases are discussed next.

In Fig. 8, we compare the two drag-coefficient signals (the CFD-simulated and the reduced-order modeled) for the fixed-cylinder case. In this case, $a_{1L} = 0.905$, $a_{1D} = 0.0016$, $a_{2D} = 0.077$, $\psi[a_{1D}, a_{1L}] = 256°$, and $\psi[a_{2D}, a_{1L}] = 335.2°$. Consequently, $\lambda_1 = -0.241$ and $\lambda_2 = -0.970$; while $q_1 = 0.908$ and $q_2 = -0.419$. The mean drag coefficient (its DC or zero-frequency component) based on the training CFD signal is 1.376. As shown in the figure, the agreement for the drag coefficient is excellent. We would like to add here that the superharmonics in the lift coefficient beyond $a_{3L}$ are negligible, which allows replacing the lift coefficient signal obtained from the CFD simulation with a third-order approximation of the form:

$$C_L \cong a_{1L} \cos(\Omega t) + a_{2L} \cos(2\Omega t + \psi[a_{2L}/a_{1L}]) \\ + a_{3L} \cos(3\Omega t + \psi[a_{3L}/a_{1L}]) \quad (8)$$

The definitions of $\psi[a_{2L}, a_{1L}]$ and $\psi[a_{3L}, a_{1L}]$ are similar to the definition of $\psi[a_{1D}, a_{1L}]$. In Fig. 9, we compare this harmonic approximation of the lift coefficient to the simulation results, and also the modeled drag coefficient that corresponds to either of the two lift signals. The difference is very small, in fact, which is again due to the insignificance of the neglected superharmonic components in the lift coefficient. Thus, Eq. 8 can be viewed as an algebraic model for the lift coefficient.

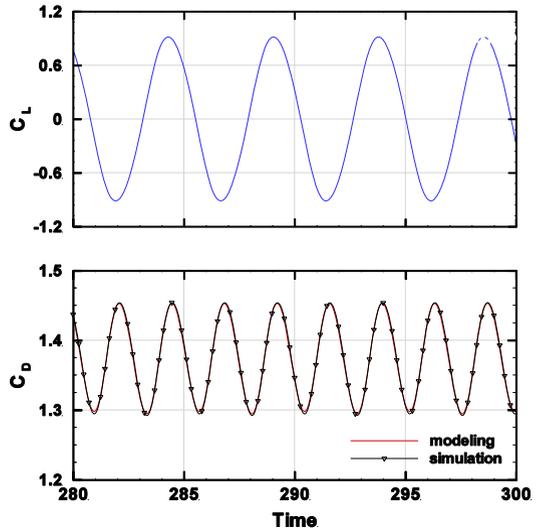

**Fig. 8** Comparison of the modeled and simulated drag coefficient for the fixed cylinder. The simulated lift coefficient is also shown

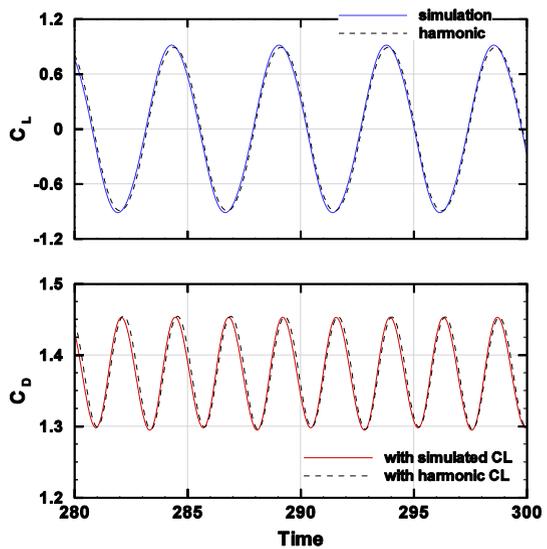

**Fig. 9** Comparison of the analytical (significant harmonics) and numerical (simulation) representation of the lift coefficient for the fixed cylinder, and the corresponding modeled drag

In Fig. 10, we compare the two drag-coefficient signals for a case where $a_{1D}$ is significant and even larger than $a_{2D}$. This case corresponds to case (b) in Fig. 6–7 (the tilt angle of the mechanical excitation line is 50° from the horizontal free stream), $a_{1L} = 0.902$, $a_{1D} = 0.16$, $a_{2D} = 0.123$, $\psi[a_{1D}, a_{1L}] = 223°$, and $\psi[a_{2D}, a_{1L}] = 352.9°$. Consequently, $\lambda_1 = -0.736$ and $\lambda_2 = -0.677$; while $q_1 = 0.992$ and $q_2 =$

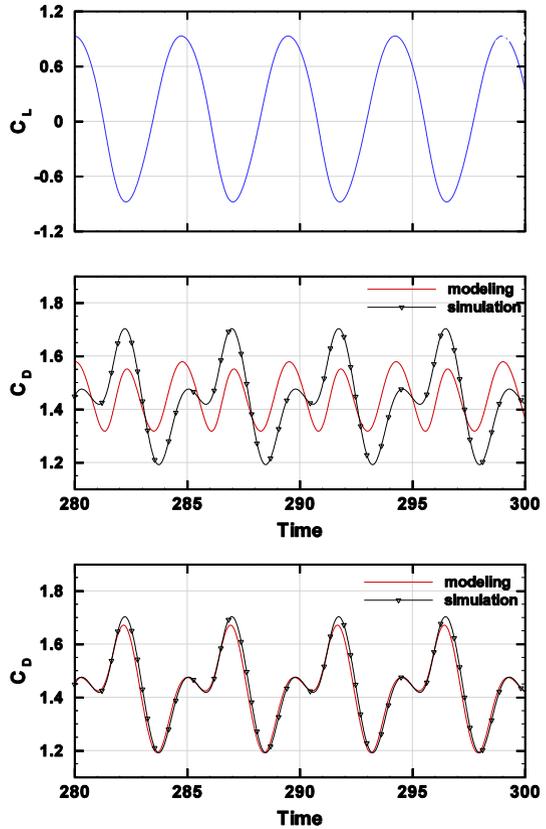

**Fig. 10** Comparison of the modeled and simulated drag coefficient for a moving cylinder with $a_{1D} > a_{2D}$, with (middle) and without (bottom) linear coupling. The simulated lift coefficient is also shown (top). Simulation data correspond to case (b) in Fig. 6–7

–0.124. The mean drag coefficient (its DC or zero-frequency component) based on the training CFD signal is 1.445. In the same figure, we also show the $C_L$ signal (as obtained from the CFD simulation). Again, the fact that the main harmonic of the lift is more than one order of magnitude larger than any of its higher harmonics (superharmonics) minimizes the modulation in the $C_L$ signal and makes it close to a pure-harmonic signal. In the same figure, we demonstrate the significance of the added linear-coupling terms by comparing the simulated drag-coefficient signal with the modeled one, but after dropping the linear terms. The disagreement is remarkable, which is a result of assuming a two-to-one frequency relationship between the lift and drag in this case, and such an assumption is not true in this case.

As the motion axis gets closer to the cross-stream, the amplitude $a_{1D}$ decreases, whereas $a_{2D}$ exhibits slight variations. For case (d) in. Figure 6–7(the tilt angle of the mechanical excitation line is 70° from the horizontal free stream), $a_{1D} = 0.0796$ and $a_{2D} = 0.124$. Although $a_{2D}$ is larger now than $a_{1D}$, the contribution of the latter is still important. In Fig. 11, we indicate this by comparing the modeled drag-coefficient signal with and without the linear-coupling terms (i.e., with and without accounting for $a_{1D}$) to the simulated signal. Again, the proposed universal drag model performs very well, whereas the original model (without the extra added terms) shows some level of discrepancy that is less than the one in Fig. 10 (due to decreased ratio of $a_{1D}/a_{2D}$), but it is still considerable.

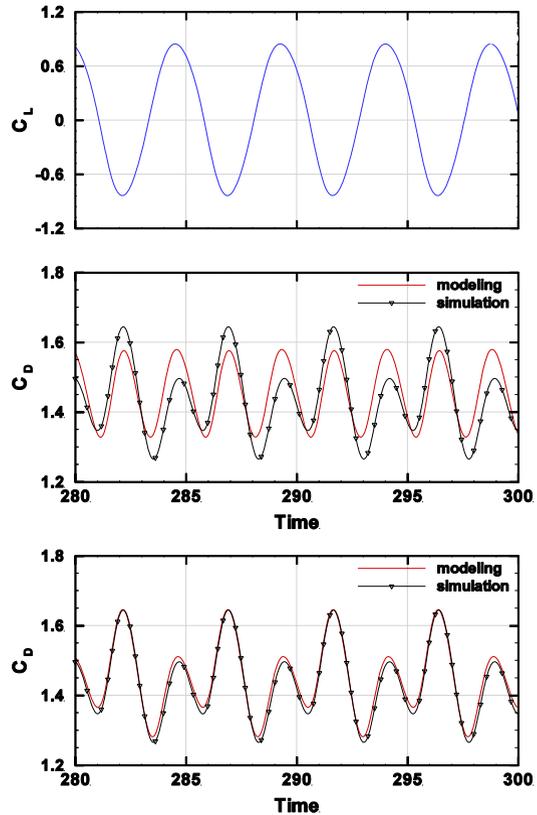

**Fig. 11** Comparison of the modeled and simulated drag coefficient for moving cylinder with $a_{1D} < a_{2D}$, with (middle) and without (bottom) linear coupling. The simulated lift coefficient is also shown (top). Simulation data correspond to case (d) in Fig. 6–7

The final case is for the cross-stream motion, which is case (f) in. Figure 6–7. This case does not represent a challenge to the proposed model since the two-to-one frequency relationship between the drag and lift is restored (as was originally the case for a fixed cylinder, before introducing any motion). This explains the big similarity between the modeled drag signal with and without the linear terms, as indicated in Fig. 12.

We would like to add that after running the simulations for several cases of interest (i.e., at different tilt angles of the cylinder's vibration axis), one can fit the various needed spectral quantities for the proposed model, such as $a_{1L}$ and $a_{1D}$, into suitable functions of the inclination angle, which can be used later to produce the modeled drag signal at any arbitrary inclination angle. This saves time remarkably, knowing that generating the modeled signal takes about three orders of magnitude less time than generating the CFD-based signal. This also allows for analyzing the tilt sensitivity [138, 139] of the flow parameters. Similarly, this concept of one-time construction of a lookup table (LUT) or database [140, 141] of model parameters can be repeated at different Reynolds numbers of interest or different oscillation amplitudes. After building such a database of parameters, the reduced modeling can be used to conveniently generate interpolated signals within a continuous parameter space.

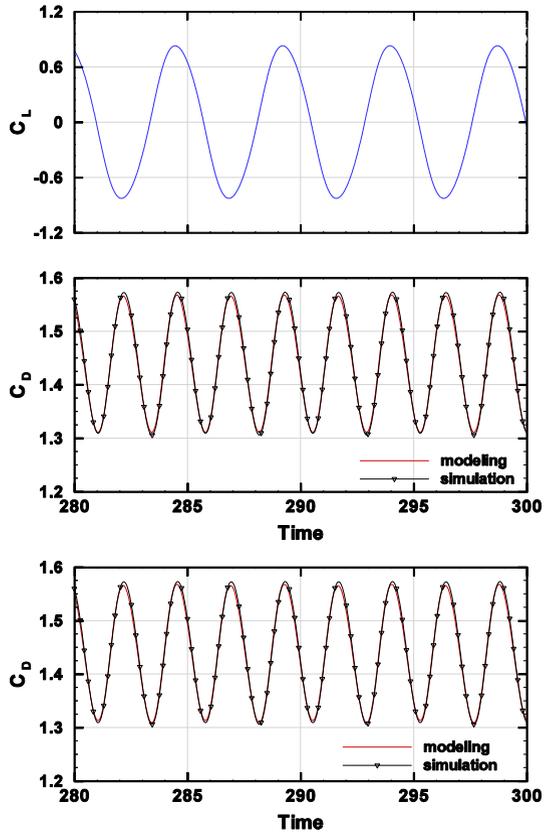

**Fig. 12** Comparison of the modeled and simulated drag coefficient for a moving cylinder in the cross-stream direction, with (middle) and without (bottom) linear coupling. The simulated lift coefficient is also shown (top). Simulation data correspond to case (f) in Fig. 6–7

We also would like to comment on the applicability of the proposed extended drag model in more complex three-dimensional flow regimes (such as in the case of a finite-length cylinder with a high Reynolds number of $10^6$, for example) and whether the two added nonlinear terms in the drag coefficient model are able to explain the drag-lift coupling in such three-dimensional flow dynamics. Due to the three-dimensional nature of such a situation, the signals of either the lift coefficient or the drag coefficient in a real problem become less regular than in the case of simpler two-dimensional flow. However, the three-dimensional flow still shows a coherent vortex street with distinguished frequency components in the lift coefficient [142, 143], in a similar manner to the two-dimensional flow problem [144, 145]. Therefore, the proposed drag model remains useful in three-dimensional problems.

## 5 Discussion

Although the current study focuses on reduced-order modeling of the drag coefficient, it can be useful to briefly discuss the reduced-order oscillator model for the lift coefficient also. The current section is devoted to achieving this aim.

Despite the availability of more than one form for modeling the lift coefficient, we consider here one simple form (van der Pol oscillator) for a fixed cylinder. Such a simple oscillator facilitates the discussion and avoids the extensive details that are involved with a more complicated modeling approach.

A van der Pol oscillator [146–148] for the lift coefficient ($C_L$) has the following form:

$$\ddot{C}_L + \omega_S C_L^2 - \mu \dot{C}_L + \alpha C_L^2 \dot{C}_L = 0 \tag{9}$$

where $\ddot{C}_L$ is the second time derivative of the lift coefficient, $\omega_S$ is the Strouhal angular frequency, $\mu$ is a coefficient for linear destabilizing damping [149, 150], $\dot{C}_L$ is the first time derivative of the lift coefficient, and $\alpha$ is a coefficient for nonlinear (cubic) stabilizing damping [151–153].

At a low value of $C_L$ (and thus $C_L^2$), the nonlinear damping term becomes weaker than the linear destabilizing damping term. Therefore, the lift coefficient tends to be amplified. On the other hand, a high value of $C_L$ (and thus $C_L^2$) causes the nonlinear damping term to be stronger than the linear destabilizing damping term. Therefore, the lift coefficient tends to decay. The opposing influences of the two damping terms lead to a balanced response that is a periodic limit cycle [154–156], where the time signal of the lift coefficient oscillates with a fixed amplitude.

## 6 Conclusions

We considered non-traditional wakes formed as a result of applied mechanical oscillatory perturbations, where the wake of a fixed cylinder is excited by introducing a prescribed vibrational motion along an inclined straight line that lies between the cross-stream direction and the stream-wise direction. This revealed new types of the relationship between the induced lift and drag forces acting on the cylinder. The typical quadratic lift-drag coupling alone fails to reproduce this relationship. We first analyzed different lift and drag signals obtained from computational fluid dynamics (CFD) simulations. The analysis techniques were in the time domain (time series analysis), the frequency domain (spectral analysis), and the lift-drag plane (limit cycle projection). Our observation that the two-loop curve of the lift-drag projection of the limit cycle gradually transforms into a single-loop curve suggested that additional linear lift-coupling terms in the drag reduced-order model (ROM) are needed, especially when the previously neglected component of the drag coefficient at the frequency of the main lift component ($a_{1D}$) increases. We added two linear terms to an existing three-term drag ROM so that it is extended and can handle non-traditional wake cases where the two-to-one frequency relationship between the drag and lift is broken. We related the model parameters to the spectral variables of the lift and drag signals. We performed several tests of the proposed five-term universal drag model, including the fixed cylinder case and the cross-stream vibration case. These two special cases exhibit a two-to-one frequency relationship between the drag and lift. We also performed tests of the proposed five-term algebraic drag model for cases where the drag has significant components at the motion frequency in addition to the typical component at twice this frequency, so both linear and quadratic couplings become important. All the comparisons showed excellent agreement between the proposed algebraic reduced-order drag model and the time-consuming CFD simulations, which require solving the Navier–Stokes equations that govern the flow field.

**Author contributions** Not applicable (this research has a single author). OM was in charge of the methodology, software, validation, formal analysis, investigation, visualization, writing of the original draft, and writing of the revised version.

**Funding** Not applicable (this research received no funding).

**Data availability** Not applicable (the data that support the findings of this study are available within the article itself).

## Declarations

**Competing interests** The authors declare no competing interests.

**Ethical approval and consent to participate** Not applicable (this research does not involve human participants, human data, human tissue, animal subjects, or environmental hazards).

**Consent for publication** Not applicable (this manuscript does not contain data from any individual person).

## References

1. Marzouk, O.; Nayfeh, A. Differential/Algebraic Wake Model Based on the Total Fluid Force and Its Direction, and the Effect of Oblique Immersed-Body Motion on 'Type-1' and 'Type-2' Lock-In. In *47th AIAA Aerospace Sciences Meeting including The New Horizons Forum and Aerospace Exposition*; AIAA [American Institute of Aeronautics and Astronautics]: Orlando, Florida, USA, 2009; p AIAA 2009–1112. https://doi.org/10.2514/6.2009-1112.
2. Gholami A, Farhang Mehr V, Basalani H (2014) Experimental study of supercavitating flow over a cylindrical body equipped with various activators. Iran J Mech Eng Trans ISME 16(2):90–109
3. Marzouk, O. A.; Nayfeh, A. H. Simulation, Analysis, and Explanation of the Lift Suppression and Break of 2:1 Force Coupling Due to in-Line Structural Vibration. In *49th AIAA/ASME/ASCE/AHS/ASC Structures, Structural Dynamics, and Materials Conference*; AIAA [American Institute of Aeronautics and Astronautics]: Schaumburg, Illinois, USA, 2008; p AIAA 2008–2309. https://doi.org/10.2514/6.2008-2309.

4. Odru, P.; Poirette, Y.; Stassen, Y.; Saint Marcoux, J. F.; Abergel, L. Composite Riser and Export Line Systems for Deep Offshore Applications; American Society of Mechanical Engineers Digital Collection, 2009; pp 147–156. https://doi.org/10.1115/OMAE2003-37237.

5. Pinto, A.; Broglia, R.; Di Mascio, A.; Campana, E. F.; Rocco, P. Numerical Investigation of the Unsteady Flow at High Reynolds Number Over a Marine Riser With Helical Strakes; American Society of Mechanical Engineers Digital Collection, 2008; pp 587–595. https://doi.org/10.1115/OMAE2006-92161.

6. Pierella F, Sætran L (2017) Wind tunnel investigation on the effect of the turbine tower on wind turbines wake symmetry. Wind Energ 20(10):1753–1769. https://doi.org/10.1002/we.2120

7. Marzouk OA (2025) Wind speed Weibull model identification in Oman, and computed normalized annual energy production (NAEP) from wind turbines based on data from weather stations. Eng Rep 7(3):e70089. https://doi.org/10.1002/eng2.70089

8. Marzouk OA (2024) Portrait of the decarbonization and renewables penetration in Oman's energy mix, motivated by Oman's National Green Hydrogen Plan. Energies 17(19):4769. https://doi.org/10.3390/en17194769

9. Harris, R. L. A Numerical Analysis of the Flow Field Surrounding a Solar Chimney Power Plant. Master of Science in Engineering, Stellenbosch : University of Stellenbosch, University of Stellenbosch, South Africa, 2004. http://hdl.handle.net/10019.1/16337 (accessed 2025–06–05).

10. Marzouk OA (2024) Energy generation intensity (EGI) of solar updraft tower (SUT) power plants relative to CSP plants and PV power plants using the new energy simulator "Aladdin." Energies 17(2):405. https://doi.org/10.3390/en17020405

11. Papailiou KO (2021) Overhead Lines. In: Papailiou KO (ed) Springer Handbook of Power Systems. Springer, Singapore, pp 611–758

12. MacIver C, Cruden A, Leithead WE, Bertinat MP (2015) Effect of wind turbine wakes on wind-induced motions in wood-pole overhead lines. Wind Energ 18(4):643–662. https://doi.org/10.1002/we.1717

13. Marzouk OA (2024) Estimated electric conductivities of thermal plasma for air-fuel combustion and oxy-fuel combustion with potassium or cesium seeding. Heliyon 10(11):e31697. https://doi.org/10.1016/j.heliyon.2024.e31697

14. Marzouk OA (2024) Subcritical and supercritical Rankine steam cycles, under elevated temperatures up to 900°C and absolute pressures up to 400 Bara. Adv Mech Eng 16(1):16878132231221065. https://doi.org/10.1177/16878132231221065

15. John AD, Gairola A, Ganju E, Gupta A (2011) Design wind loads on reinforced concrete chimney – an experimental case study. Procedia Eng 14:1252–1257. https://doi.org/10.1016/j.proeng.2011.07.157

16. Sun Y, Li Z, Sun X, Su N, Peng S (2020) Interference effects between two tall chimneys on wind loads and dynamic responses. J Wind Eng Ind Aerodyn 206:104227. https://doi.org/10.1016/j.jweia.2020.104227

17. Marzouk OA (2024) Expectations for the role of hydrogen and its derivatives in different sectors through analysis of the four energy scenarios: IEA-STEPS, IEA-NZE, IRENA-PES, and IRENA-1.5°C. Energies 17(3):646. https://doi.org/10.3390/en17030646

18. Easter, N.; Gao, M.; Krishnamurthy, R. M.; Kompally, S.; Coates, T. Implication of Fractographic Analysis of the Crack in an Above Ground Pipeline – Potential Role of Hydrogen; OnePetro, 2022.

19. Marzouk OA (2023) 2030 ambitions for hydrogen, clean hydrogen, and green hydrogen. Eng Proc 56(1):14. https://doi.org/10.3390/ASEC2023-15497

20. Marzouk OA (2023) Zero carbon ready metrics for a single-family home in the Sultanate of Oman based on EDGE certification system for green buildings. Sustainability 15(18):13856. https://doi.org/10.3390/su151813856

21. *Buckling of Thin Metal Shells*; Teng, J. G., Rotter, J. M., Eds.; CRC Press: London, UK, 2006. https://doi.org/10.1201/9781482295078.

22. Marzouk OA (2022) Land-use competitiveness of photovoltaic and concentrated solar power technologies near the Tropic of Cancer. Sol Energy 243:103–119. https://doi.org/10.1016/j.solener.2022.07.051

23. Marzouk OA (2025) Solar heat for industrial processes (SHIP): an overview of its categories and a review of its recent progress. Solar 5(4):46. https://doi.org/10.3390/solar5040046

24. Al Manthari MS, Azeez C, Sankar M, Pushpa BV (2024) Numerical study of laminar flow and vortex-induced vibration on cylinder subjects to free and forced oscillation at low Reynolds numbers. Fluids 9(8):175. https://doi.org/10.3390/fluids9080175
25. Lamont PJ, Hunt BL (1973) Out-of-plane force on a circular cylinder at large angles of inclination to a uniform stream. Aeronaut J 77(745):41–45. https://doi.org/10.1017/S0001924000040136
26. Civrais CHB, White C, Steijl R (2024) Influence of Anharmonic oscillator model for flows over a cylindrical body. AIP Conf Proc 2996(1):080008. https://doi.org/10.1063/5.0187445
27. Marzouk OA (2025) Reduced-order modeling (ROM) of a segmented plug-flow reactor (PFR) for hydrogen separation in integrated gasification combined cycles (IGCC). Processes 13(5):1455. https://doi.org/10.3390/pr13051455
28. Zhang W, Kou J, Wang Z (2016) Nonlinear aerodynamic reduced-order model for limit-cycle oscillation and flutter. AIAA J 54(10):3304–3311. https://doi.org/10.2514/1.J054951
29. Marzouk, O. A.; Nayfeh, A. H. Detailed Characteristics of the Resonating and Non-Resonating Flows Past a Moving Cylinder. In *49th AIAA/ASME/ASCE/AHS/ASC Structures, Structural Dynamics, and Materials Conference*; AIAA [American Institute of Aeronautics and Astronautics]: Schaumburg, Illinois, USA, 2008; p AIAA 2008–2311. https://doi.org/10.2514/6.2008-2311
30. Bishop RED, Hassan AY (1997) The lift and drag forces on a circular cylinder in a flowing fluid. Proc R Soc Lond Ser A Math Phys Sci 277(1368):32–50. https://doi.org/10.1098/rspa.1964.0004
31. Marzouk OA (2022) Urban air mobility and flying cars: overview, examples, prospects, drawbacks, and solutions. Open Eng 12(1):662–679. https://doi.org/10.1515/eng-2022-0379
32. Vickery BJ (1966) Fluctuating lift and drag on a long cylinder of square cross-section in a smooth and in a turbulent stream. J Fluid Mech 25(3):481–494. https://doi.org/10.1017/S002211206600020X
33. Marzouk OA (2008) A two-step computational aeroacoustics method applied to high-speed flows. Noise Control Eng J 56(5):396. https://doi.org/10.3397/1.2978229
34. Marzouk, O. A. Evolutionary Computing Applied to Design Optimization. In *ASME 2007 International Design Engineering Technical Conferences and Computers and Information in Engineering Conference (IDETC-CIE 2007), (4–7 September 2007)*; ASME [American Society of Mechanical Engineers]: Las Vegas, Nevada, USA, 2009; Vol. 2, pp 995–1003. https://doi.org/10.1115/DETC2007-35502
35. Marzouk OA (2025) Coupled differential-algebraic equations framework for modeling six-degree-of-freedom flight dynamics of asymmetric fixed-wing aircraft. Int J Adv Appl Sci 12(1):30–51. https://doi.org/10.21833/ijaas.2025.01.004
36. Massaro D, Karp M, Jansson N, Markidis S, Schlatter P (2024) Direct numerical simulation of the turbulent flow around a Flettner rotor. Sci Rep 14(1):3004. https://doi.org/10.1038/s41598-024-53194-x
37. Esmaeili M, Fakhri Vayqan H, Rabiee AH (2025) Impact of temperature-induced buoyancy on the 2DOF-VIV of a heated/cooled cylinder. Arab J Sci Eng 50(4):2807–2822. https://doi.org/10.1007/s13369-024-09262-5
38. Nagata T, Shigeta T, Kasai M, Nonomura T (2025) Schlieren visualization and drag measurement on compressible flow over a circular cylinder at Reynolds number of $\mathcal{O}(10^2)$. Exp Fluids 66(5):91. https://doi.org/10.1007/s00348-025-04010-3
39. Sridhar MK, Kang C-K, Landrum DB, Aono H, Mathis SL, Lee T (2021) Effects of flight altitude on the lift generation of Monarch butterflies: from sea level to overwintering mountain. Bioinspir Biomim 16(3):034002. https://doi.org/10.1088/1748-3190/abe108
40. Marzouk OA (2011) One-way and two-way couplings of CFD and structural models and application to the wake-body interaction. Appl Math Model 35(3):1036–1053. https://doi.org/10.1016/j.apm.2010.07.049
41. Marzouk, O. A.; Nayfeh, A. H. New Wake Models With Capability of Capturing Nonlinear Physics. In *ASME 2008 27th International Conference on Offshore Mechanics and Arctic Engineering (OMAE 2008)*; ASME [American Society of Mechanical Engineers]: Estoril, Portugal, 2009; pp 901–912. https://doi.org/10.1115/OMAE2008-57714
42. Hartlen RT, Currie IG (1970) Lift-oscillator model of vortex-induced vibration. J Eng Mech Div 96(5):577–591. https://doi.org/10.1061/JMCEA3.0001276
43. Leontini JS, Jacono DL, Thompson MC (2011) A numerical study of an inline oscillating cylinder in a free stream. J Fluid Mech 688:551–568. https://doi.org/10.1017/jfm.2011.403
44. Marzouk OA, Nayfeh AH (2009) Reduction of the loads on a cylinder undergoing harmonic in-line motion. Phys Fluids 21(8):083103. https://doi.org/10.1063/1.3210774

45. Nobari MRH, Naderan H (2006) A numerical study of flow past a cylinder with cross flow and inline oscillation. Comput Fluids 35(4):393–415. https://doi.org/10.1016/j.compfluid.2005.02.004
46. Marzouk OA, Nayfeh AH (2010) Characterization of the flow over a cylinder moving harmonically in the cross-flow direction. Int J Non-Linear Mech 45(8):821–833. https://doi.org/10.1016/j.ijnonlinmec.2010.06.004
47. D'Urso B, Van Handel R, Odom B, Hanneke D, Gabrielse G (2005) Single-particle self-excited oscillator. Phys Rev Lett 94(11):113002. https://doi.org/10.1103/PhysRevLett.94.113002
48. Marzouk OA (2025) InvSim algorithm for pre-computing airplane flight controls in limited-range autonomous missions, and demonstration via double-roll maneuver of Mirage III fighters. Sci Rep 15:23382. https://doi.org/10.1038/s41598-025-07639-6
49. Kappauf J, Hetzler H (2021) On a hybrid approximation concept for self-excited periodic oscillations of large-scale dynamical systems. PAMM 21(1):e202100143. https://doi.org/10.1002/pamm.202100143
50. Skop RA, Griffin OM (1973) A model for the vortex-excited resonant response of bluff cylinders. J Sound Vib 27(2):225–233. https://doi.org/10.1016/0022-460X(73)90063-1
51. Blevins, R. D. Flow Induced Vibration of Bluff Structures. PhD in Engineering Mechanics, California Institute of Technology (Caltech), Pasadena, California, USA, 1974. https://www.proquest.com/openview/bed7d14a974e2b448029fb25620d77f3 (accessed 2025–06–03).
52. Di Silvio G, Angrilli F, Zanardo A (1975) Fluidelastic vibrations: mathematical model and experimental result. Meccanica 10(4):269–279. https://doi.org/10.1007/BF02133219
53. Landl R (1975) A mathematical model for vortex-excited vibrations of bluff bodies. J Sound Vib 42(2):219–234. https://doi.org/10.1016/0022-460X(75)90217-5
54. Goswami I, Scanlan RH, Jones NP (1993) Vortex-induced vibration of circular cylinders. II: new model. J Eng Mech 119(11):2288–2302. https://doi.org/10.1061/(ASCE)0733-9399(1993)119:11(2288)
55. Skop RA, Balasubramanian S (1997) A new twist on an old model for vortex-excited vibrations. J Fluids Struct 11(4):395–412. https://doi.org/10.1006/jfls.1997.0085
56. Krenk S, Nielsen SRK (1999) Energy balanced double oscillator model for vortex-induced vibrations. J Eng Mech 125(3):263–271. https://doi.org/10.1061/(ASCE)0733-9399(1999)125:3(263)
57. Mureithi NW, Goda S, Kanki H, Nakamura T (2001) A nonlinear dynamics analysis of vortex-structure interaction models. J Press Vessel Technol 123(4):475–479. https://doi.org/10.1115/1.1403023
58. Facchinetti ML, de Langre E, Biolley F (2004) Coupling of structure and wake oscillators in vortex-induced vibrations. J Fluids Struct 19(2):123–140. https://doi.org/10.1016/j.jfluidstructs.2003.12.004
59. Marzouk OA (2009) A nonlinear ODE system for the unsteady hydrodynamic force - a new approach. World Acad Sci Eng Tech 39:948–962.
60. Kim W-J, Perkins NC (2002) Two-dimensional vortex-induced vibration of cable suspensions. J Fluids Struct 16(2):229–245. https://doi.org/10.1006/jfls.2001.0418
61. Ren, C.; Cheng, L.; Tong, F. Hydrodynamic Force of an Obliquely Oscillating Cylinder in Steady Flow; American Society of Mechanical Engineers Digital Collection, 2024. https://doi.org/10.1115/OMAE2024-126997.
62. Marzouk, O. A.; Nayfeh, A. H. Fluid Forces and Structure-Induced Damping of Obliquely-Oscillating Offshore Structures. In *The Eighteenth International Offshore and Polar Engineering Conference (ISOPE-2008)*; ISOPE [International Society of Offshore and Polar Engineers]: Vancouver, British Columbia, Canada, 2008; Vol. 3, pp 460–468.
63. Wang XQ, So RMC, Chan KT (2003) A non-linear fluid force model for vortex-induced vibration of an elastic cylinder. J Sound Vib 260(2):287–305. https://doi.org/10.1016/S0022-460X(02)00945-8
64. Awrejcewicz J (2017). Resonance. https://doi.org/10.5772/intechopen.68248
65. Thomas, F.; Chaney, R.; Tseng, R. Resonance. In *The Physics of Destructive Earthquakes*; Morgan & Claypool Publishers: San Rafael, California, USA, 2018; Vol. The Physics of Destructive Earthquakes, pp 8–11. https://doi.org/10.1088/978-1-64327-078-4ch8.
66. Mignolet, M.; Red-Horse, J. ARMAX Identification of Vibrating Structures - Model and Model Order Estimation. In *35th Structures, Structural Dynamics, and Materials Conference*; AIAA [American Institute of Aeronautics and Astronautics]: Hilton Head, South Carolina, USA, 1994; p AIAA-94–1525-CP. https://doi.org/10.2514/6.1994-1525.
67. Chen-xu, N.; Jie-sheng, W. Auto Regressive Moving Average (ARMA) Prediction Method of Bank Cash Flow Time Series. In *2015 34th Chinese Control Conference (CCC)*; 2015; pp 4928–4933. https://doi.org/10.1109/ChiCC.2015.7260405.

68. Huang L (2015) Auto regressive moving average (ARMA) modeling method for gyro random noise using a robust Kalman filter. Sensors 15(10):25277–25286. https://doi.org/10.3390/s151025277
69. Liao Y, Tang H, Xie L (2025) A deep modal model for reconstructing the VIV response of a flexible cylinder with sparse sensing data. Ocean Eng 326:120871. https://doi.org/10.1016/j.oceaneng.2025.120871
70. Marzouk OA (2010) Characteristics of the flow-induced vibration and forces with 1- and 2-DOF vibrations and limiting solid-to-fluid density ratios. J Vib Acoust 132(4):041013. https://doi.org/10.1115/1.4001503
71. Srinivasan K, Prethiv Kumar R, Nallayarasu S, Liu Y (2025) Frequency-domain analysis of vortex-induced vibration of flexible cantilever cylinder with various aspect ratios. Ocean Eng 320:120204. https://doi.org/10.1016/j.oceaneng.2024.120204
72. Wallar, B. L.; Kimber, M. L. Numerical Investigation of Force Coefficient Data for Multiple Flow Past Cylinder Configurations. In *10th Thermal and Fluids Engineering Conference (TFEC)*; Begel House Inc.: Washington, D.C., USA, 2025; pp 469–478. https://doi.org/10.1615/TFEC2025.bio.055994.
73. Marzouk OA (2018) Radiant heat transfer in nitrogen-free combustion environments. Int J Nonlinear Sci Numer Simul 19(2):175–188. https://doi.org/10.1515/ijnsns-2017-0106
74. Amer MN, Abuelyamen A, Parezanović VB, Alkaabi AK, Alameri SA, Afgan I (2025) A comprehensive review, CFD and ML analysis of flow around tandem circular cylinders at sub-critical Reynolds numbers. J Wind Eng Ind Aerodyn 257:105998. https://doi.org/10.1016/j.jweia.2024.105998
75. Marzouk OA (2017) Performance analysis of shell-and-tube dehydrogenation module: dehydrogenation module. Int J Energy Res 41(4):604–610. https://doi.org/10.1002/er.3637
76. Grimm V, Heinlein A, Klawonn A (2025) Learning the solution operator of two-dimensional incompressible Navier-Stokes equations using physics-aware convolutional neural networks. J Comput Phys 535:114027. https://doi.org/10.1016/j.jcp.2025.114027
77. Hoffmann, K. A.; Chiang, S. T. *Computational Fluid Dynamics - Volume 1*, 4. ed., 2. print.; Engineering Education System: Wichita, Kansas, USA, 2004.
78. Zhang Y, Rabczuk T, Lin J, Lu J, Chen CS (2024) Numerical simulations of two-dimensional incompressible Navier-Stokes equations by the backward substitution projection method. Appl Math Comput 466:128472. https://doi.org/10.1016/j.amc.2023.128472
79. Marzouk OA (2023) Adiabatic flame temperatures for Oxy-methane, Oxy-hydrogen, air-methane, and air-hydrogen stoichiometric combustion using the NASA CEARUN Tool, GRI-Mech 3.0 reaction mechanism, and Cantera python package. Eng Tech Appl Sci Res 13(4):11437–11444. https://doi.org/10.48084/etasr.6132
80. Marzouk OA (2025) Dataset of total emissivity for CO2, H2O, and H2O-CO2 mixtures; over a temperature range of 300–2900 K and a pressure-pathlength range of 0.01–50 Atm.m. Data Brief 59:111428. https://doi.org/10.1016/j.dib.2025.111428
81. Marzouk OA (2023) Temperature-dependent functions of the electron–neutral momentum transfer collision cross sections of selected combustion plasma species. Appl Sci 13(20):11282. https://doi.org/10.3390/app132011282
82. Marzouk OA (2025) Power density and thermochemical properties of hydrogen magnetohydrodynamic (H2MHD) generators at different pressures, seed types, seed levels, and oxidizers. Hydrogen 6(2):31. https://doi.org/10.3390/hydrogen6020031
83. Dondapati RS, Rao VV (2012) CFD analysis of cable-in-conduit conductors (CICC) for fusion grade magnets. IEEE Trans Appl Supercond 22(3):4703105–4703105. https://doi.org/10.1109/TASC.2012.2185025
84. Hah C (1986) A numerical modeling of endwall and tip-clearance flow of an isolated compressor rotor. J Eng Gas Turbines Power 108(1):15–21. https://doi.org/10.1115/1.3239863
85. Wu J, Liu Y, Zhang D (2025) Numerical simulation of flow around a transversely oscillating square cylinder at different frequencies. Phys Fluids 37(2):025213. https://doi.org/10.1063/5.0256326
86. Marzouk OA (2009) Direct numerical simulations of the flow past a cylinder moving with sinusoidal and nonsinusoidal profiles. J Fluids Eng 131(12):121201. https://doi.org/10.1115/1.4000406
87. Marzouk OA (2020) The sod gasdynamics problem as a tool for benchmarking face flux construction in the finite volume method. Sci Afr 10:e00573. https://doi.org/10.1016/j.sciaf.2020.e00573
88. Marzouk OA (2025) OpenFOAM computational fluid dynamics (CFD) solver for magnetohydrodynamic open cycles, applied to the Sakhalin pulsed magnetohydrodynamic generator (PMHDG). Discover Appl Sci 7(10):1108. https://doi.org/10.1007/s42452-025-07744-1

89. Marzouk, O. A. Directivity and Noise Propagation for Supersonic Free Jets. In *46th AIAA Aerospace Sciences Meeting and Exhibit*; AIAA [American Institute of Aeronautics and Astronautics]: Reno, Nevada, USA, 2008; p AIAA 2008–23. https://doi.org/10.2514/6.2008-23.
90. Vos, R.; Vaessen, F. A New Compressibility Correction Method to Predict Aerodynamic Interaction between Lifting Surfaces. In *2013 Aviation Technology, Integration, and Operations Conference*; AIAA [American Institute of Aeronautics and Astronautics]: Los Angeles, California, USA, 2013; p AIAA 2013–4299. https://doi.org/10.2514/6.2013-4299.
91. Sultanian, B. K. Fluid Mechanics and Turbomachinery: Problems and Solutions; CRC Press: Boca Raton, Florida, USA, 2021. https://doi.org/10.1201/9781003053996.
92. Marzouk OA (2025) Detailed derivation of the scalar explicit expressions governing the electric field, current density, and volumetric power density in the four types of linear divergent MHD channels under a unidirectional applied magnetic field. Contemp Math 6(4):4060–4100. https://doi.org/10.37256/cm.6420256918
93. Crowdy DG (2006) Analytical solutions for uniform potential flow past multiple cylinders. Eur J Mech B Fluids 25(4):459–470. https://doi.org/10.1016/j.euromechflu.2005.11.005
94. Marzouk OA (2023) Detailed and simplified plasma models in combined-cycle magnetohydrodynamic power systems. Int J Adv Appl Sci 10(11):96–108. https://doi.org/10.21833/ijaas.2023.11.013
95. Pindado S, Meseguer J (2003) Wind tunnel study on the influence of different parapets on the roof pressure distribution of low-rise buildings. J Wind Eng Ind Aerodyn 91(9):1133–1139. https://doi.org/10.1016/S0167-6105(03)00055-2
96. Marzouk OA (2021) Assessment of global warming in Al Buraimi, sultanate of oman based on statistical analysis of NASA power data over 39 years, and testing the reliability of NASA POWER against meteorological measurements. Heliyon 7(3):e06625. https://doi.org/10.1016/j.heliyon.2021.e06625
97. Batten WMJ, Bahaj AS, Molland AF, Chaplin JR (2006) Hydrodynamics of marine current turbines. Renewable Energy 31(2):249–256. https://doi.org/10.1016/j.renene.2005.08.020
98. Alexander DE (2017) Fluid Biomechanics. In Nature's Machines; Elsevier, Amsterdam, pp 51–97
99. Kida S, Miura H (1998) Swirl condition in low-pressure vortices. J Phys Soc Jpn 67(7):2166–2169. https://doi.org/10.1143/JPSJ.67.2166
100. Kida, S.; Goto, S.; Makihara, T. Low-Pressure Vortex Analysis in Turbulence: Life, Structure, and Dynamical Role of Vortices. In *Tubes, Sheets and Singularities in Fluid Dynamics*; Bajer, K., Moffatt, H. K., Eds.; Springer Netherlands: Dordrecht, 2002; pp 181–190. https://doi.org/10.1007/0-306-48420-X_25.
101. Marzouk OA (2024) Hydrogen utilization as a plasma source for magnetohydrodynamic direct power extraction (MHD-DPE). IEEE Access 12:167088–167107. https://doi.org/10.1109/ACCESS.2024.3496796
102. Zdravkovich M (1988) Conceptual overview of laminar and turbulent flows past smooth and rough circular cylinders. Wind Engineers, JAWE 1988(37):93–102. https://doi.org/10.5359/jawe.1988.37_93
103. Guohui H, Dejun S, Xieyuan Y (1997) On the topological bifurcation of flows around a rotating circular cylinder. Acta Mech Sinica 13(3):203–209. https://doi.org/10.1007/BF02487702
104. Wei, J.; Wu, J. Flow-Induced Reconfiguration of and Force on Elastic Cantilevers. In *Proceedings of the IUTAM Symposium on Turbulent Structure and Particles-Turbulence Interaction*; Zheng, X., Balachandar, S., Eds.; Springer Nature Switzerland: Cham, Switzerland, 2024; pp 229–249. https://doi.org/10.1007/978-3-031-47258-9_15.
105. Lehmkuhl O, Rodríguez I, Borrell R, Chiva J, Oliva A (2014) Unsteady forces on a circular cylinder at critical reynolds numbers. Phys Fluids 26(12):125110. https://doi.org/10.1063/1.4904415
106. Gowen, F. E.; Perkins, E. W. *DRAG OF CIRCULAR CYLINDERS FOR A WIDE RANGE OF REYNOLDS NUMBERS AND MACH NUMBERS*; RESEARCH MEMORANDUM RM A52C20; NACA [United States National Advisory Committee for Aeronautics]: Washington, D.C, USA, 1952; pp 1–28. https://ntrs.nasa.gov/api/citations/19930087134/downloads/19930087134.pdf (accessed 2025-11-06).
107. Kološ I, Michalcová V, Lausová L (2021) Numerical analysis of flow around a cylinder in critical and subcritical regime. Sustainability 13(4):2048. https://doi.org/10.3390/su13042048
108. Warren BA, LaCasce JH, Robbins PE (1996) On the obscurantist physics of "Form Drag" in theorizing about the circumpolar current. J Phys Oceanogr 26(10):2297–2301
109. Emeis S (1990) Pressure drag of obstacles in the Atmospheric Boundary Layer. J Appl Meteorol Climatol 29(6):461–476. https://doi.org/10.1175/1520-0450(1990)029%3c0461:PDOOIT%3e2.0.CO;2

110. Wang J, Yang Q, Huang B, Ouyang Y (2024) Adaptive ANCF method described by arbitrary Lagrange-Euler formulation with application in variable-length underwater tethered systems moving in limited spaces. Ocean Eng 297:117059. https://doi.org/10.1016/j.oceaneng.2024.117059
111. Marzouk OA (2011) Flow control using bifrequency motion. Theor Comput Fluid Dyn 25(6):381–405. https://doi.org/10.1007/s00162-010-0206-6
112. Li B, Ma S, Qiu W (2025) Optimal convergence of the arbitrary Lagrangian-Eulerian interface tracking method for two-phase navier-stokes flow without surface tension. IMA J Numer Anal. https://doi.org/10.1093/imanum/draf003
113. Marzouk OA (2010) Contrasting the Cartesian and polar forms of the shedding-induced force vector in response to 12 subharmonic and superharmonic mechanical excitations. Fluid Dyn Res 42(3):035507. https://doi.org/10.1088/0169-5983/42/3/035507
114. Marzouk OA (2024) Benchmarks for the Omani higher education students-faculty ratio (SFR) based on World Bank data, QS rankings, and THE rankings. Cogent Educ 11(1):2317117. https://doi.org/10.1080/2331186X.2024.2317117
115. Prasad, S.; G, B. Benchmarking in CFD; American Society of Mechanical Engineers Digital Collection, 2009; pp 1221–1227. https://doi.org/10.1115/HT-FED2004-56746.
116. Marzouk OA (2025) Benchmarking retention, progression, and graduation rates in undergraduate higher education across different time windows. Cogent Educ 12(1):2498170. https://doi.org/10.1080/2331186X.2025.2498170
117. Blackburn HM, Henderson RD (1999) A study of two-dimensional flow past an oscillating cylinder. J Fluid Mech 385:255–286. https://doi.org/10.1017/S0022112099004309
118. Bhrawy AH, Baleanu D (2013) A spectral Legendre–Gauss–Lobatto collocation method for a space-fractional advection diffusion equations with variable coefficients. Rep Math Phys 72(2):219–233. https://doi.org/10.1016/S0034-4877(14)60015-X
119. Roshko, A. *On the Development of Turbulent Wakes from Vortex Streets*; Technical Note NACA-TN-2913; NACA [United States National Advisory Committee for Aeronautics]: Washington, D.C., USA, 1953; pp 1–77. https://ntrs.nasa.gov/api/citations/19930083620/downloads/19930083620.pdf (accessed 2025–06–03).
120. Wen C-Y, Yeh C-L, Wang M-J, Lin C-Y (2004) On the drag of two-dimensional flow about a circular cylinder. Phys Fluids 16(10):3828–3831. https://doi.org/10.1063/1.1789071
121. Marzouk, O. A.; Nayfeh, A. H. Loads on a Harmonically Oscillating Cylinder. In *ASME 2007 International Design Engineering Technical Conferences and Computers and Information in Engineering Conference (IDETC-CIE 2007)*; ASME [American Society of Mechanical Engineers]: Las Vegas, Nevada, USA, 2009; pp 1755–1774. https://doi.org/10.1115/DETC2007-35562.
122. Marzouk, O. A.; Nayfeh, A. H. A Study of the Forces on an Oscillating Cylinder. In *ASME 2007 26th International Conference on Offshore Mechanics and Arctic Engineering (OMAE 2007)*; ASME [American Society of Mechanical Engineers]: San Diego, California, USA, 2009; pp 741–752. https://doi.org/10.1115/OMAE2007-29163.
123. Qin, L. Development of Reduced-Order Models for Lift and Drag on Oscillating Cylinders with Higher-Order Spectral Moments. PhD in Engineering Mechanics, Virginia Polytechnic Institute and State University (Virginia Tech), Blacksburg, Virginia, USA, 2004. http://hdl.handle.net/10919/29542 (accessed 2025–06–03).
124. Modarres-Sadeghi Y (2021) Flow around a fixed cylinder. In: Modarres-Sadeghi Y (ed) Introduction to fluid-structure interactions. Springer International Publishing, Cham, Switzerland, pp 5–22
125. Bookey, P.; Wyckham, C.; Smits, A.; Martin, P. New Experimental Data of STBLI at DNS/LES Accessible Reynolds Numbers. In *43rd AIAA Aerospace Sciences Meeting and Exhibit*; AIAA [American Institute of Aeronautics and Astronautics]: Reno, Nevada, USA, 2005; p AIAA 2005–309. https://doi.org/10.2514/6.2005-309.
126. Wang Z, Yeo KS, Khoo BC (2006) DNS of low Reynolds number turbulent flows in dimpled channels. J Turbul 7:N37. https://doi.org/10.1080/14685240600595735
127. Marzouk OA, Huckaby ED (2010) A comparative study of eight finite-rate chemistry kinetics for Co/$H_2$ combustion. Eng Appl Comput Fluid Mech 4(3):331–356. https://doi.org/10.1080/19942060.2010.11015322
128. Catalano P, Tognaccini R (2010) Turbulence modeling for low-Reynolds-number flows. AIAA J 48(8):1673–1685. https://doi.org/10.2514/1.J050067

129. Marzouk OA (2025) Technical review of radiative-property modeling approaches for gray and nongray radiation, and a recommended optimized WSGGM for CO2/H2O-enriched gases. Results Eng 25:103923. https://doi.org/10.1016/j.rineng.2025.103923
130. Ou Z, Cheng L (2003) Numerical simulation of three-dimensional flow around a rectangular cylinder. In: Armfield SW, Morgan P, Srinivas K (eds) Computational fluid dynamics 2002. Springer, Berlin, Heidelberg, pp 775–776
131. Marzouk OA, Huckaby ED (2010) Simulation of a swirling gas-particle flow using different k-epsilon models and particle-parcel relationships. Eng Lett 18(1)
132. Marzouk OA, Nayfeh AH (2009) Control of ship roll using passive and active anti-roll tanks. Ocean Eng 36(9):661–671. https://doi.org/10.1016/j.oceaneng.2009.03.005
133. Marzouk, O. A.; Nayfeh, A. H. Mitigation of Ship Motion Using Passive and Active Anti-Roll Tanks. In *ASME 2007 International Design Engineering Technical Conferences and Computers and Information in Engineering Conference (IDETC-CIE 2007)*; ASME [American Society of Mechanical Engineers]: Las Vegas, Nevada, USA, 2009; pp 215–229. https://doi.org/10.1115/DETC2007-35571
134. Franke R, Rodi W, Schönung B (1990) Numerical calculation of laminar vortex-shedding flow past cylinders. J Wind Eng Ind Aerodyn 35:237–257. https://doi.org/10.1016/0167-6105(90)90219-3
135. Dong H, Fang S, Du X (2025) Flow around square, rounded, and round-convex cylinders at Reynolds numbers 20 to 22,000. Comput Fluids 300:106771. https://doi.org/10.1016/j.compfluid.2025.106771
136. Chen B, Zhu Y, Hu J, Principe JC (2013) System parameter identification: information criteria and algorithms, 1st edn. Elsevier, London
137. Marzouk OA (2022) Compilation of smart cities attributes and quantitative identification of mismatch in rankings. J Eng 2022:5981551. https://doi.org/10.1155/2022/5981551
138. Mariotti A, Galletti C, Brunazzi E, Salvetti MV (2022) Mixing sensitivity to the inclination of the lateral walls in a T-mixer. Chem Eng Process - Process Intensif 170:108699. https://doi.org/10.1016/j.cep.2021.108699
139. Marzouk OA (2022) Tilt sensitivity for a scalable one-hectare photovoltaic power plant composed of parallel racks in Muscat. Cogent Eng 9(1):2029243. https://doi.org/10.1080/23311916.2022.2029243
140. Marzouk OA (2021) Lookup tables for power generation performance of photovoltaic systems covering 40 geographic locations (Wilayats) in the Sultanate of Oman, with and without solar tracking, and general perspectives about solar irradiation. Sustainability 13(23):13209. https://doi.org/10.3390/su132313209
141. Yokoyama C, Hamada A, Ikuta Y, Shige S, Yamaji M, Tsuji H, Kubota T, Takayabu YN (2025) Spectral latent heating retrieval for the midlatitudes using GPM DPR. Part I: construction of lookup tables. J Appl Meteorol Climatol. https://doi.org/10.1175/JAMC-D-23-0217.1
142. Peppa S, Kaiktsis L, Frouzakis CE, Triantafyllou GS (2021) Computational study of three-dimensional flow past an oscillating cylinder following a figure eight trajectory. Fluids 6(3):107. https://doi.org/10.3390/fluids6030107
143. Bhagat KC, Soren SK, Chaudhary SK (2016) Experimental and numerical analysis of different aerodynamic properties of circular cylinder. Int Res J Eng Techn 3(9):1112–1117
144. Prosser, D.; Smith, M. Aerodynamics of Finite Cylinders in Quasi-Steady Flow. In *53rd AIAA Aerospace Sciences Meeting*; AIAA SciTech Forum; AIAA [American Institute of Aeronautics and Astronautics]: Kissimmee, Florida, USA, 2015; p AIAA 2015–1931. https://doi.org/10.2514/6.2015-1931
145. Chabert d'Hières G, Davies PA, Didelle H (1990) Experimental studies of lift and drag forces upon cylindrical obstacles in homogeneous, rapidly rotating fluids. Dyn Atmos Oceans 15(1):87–116. https://doi.org/10.1016/0377-0265(90)90005-G
146. Cveticanin L (2013) On the Van Der Pol oscillator: an overview. Appl Mech Mater 430:3–13. https://doi.org/10.4028/www.scientific.net/AMM.430.3
147. Mettin R, Parlitz U, Lauterborn W (1993) Bifurcation structure of the driven van Der Pol oscillator. Int J Bifurcation Chaos 03(06):1529–1555. https://doi.org/10.1142/S0218127493001203
148. Kovacic I (2011) On the motion of a generalized van Der Pol oscillator. Commun Nonlinear Sci Numer Simul 16(3):1640–1649. https://doi.org/10.1016/j.cnsns.2010.06.016
149. Sarkar A, Mondal J, Chatterjee S (2020) Controlling self-excited vibration using positive position feedback with time-delay. J Braz Soc Mech Sci Eng 42(9):464. https://doi.org/10.1007/s40430-020-02544-7

150. de Langre E (2006) Frequency lock-in is caused by coupled-mode flutter. J Fluids Struct 22(6):783–791. https://doi.org/10.1016/j.jfluidstructs.2006.04.008
151. Zainal AA, Ling CCD, Faisal SY (2011) Bifurcation of rupture path by nonlinear damping force. Appl Math Mech-Engl Ed 32(3):285–292. https://doi.org/10.1007/s10483-011-1414-9
152. Lupi F, Niemann H-J, Höffer R (2018) Aerodynamic damping model in vortex-induced vibrations for wind engineering applications. J Wind Eng Ind Aerodyn 174:281–295. https://doi.org/10.1016/j.jweia.2018.01.006
153. Zaitsev S, Shtempluck O, Buks E, Gottlieb O (2012) Nonlinear damping in a micromechanical oscillator. Nonlinear Dyn 67(1):859–883. https://doi.org/10.1007/s11071-011-0031-5
154. Robinett I, Rush D, Wilson DG (2008) What is a limit cycle? Int J Control 81(12):1886–1900. https://doi.org/10.1080/00207170801927163
155. Zhu X-M, Yin L, Ao P (2006) Limit cycle and conserved dynamics. Int J Mod Phys B 20(07):817–827. https://doi.org/10.1142/S0217979206033607
156. Thothadri M, Moon FC (2005) Nonlinear system identification of systems with periodic limit-cycle response. Nonlinear Dyn 39(1):63–77. https://doi.org/10.1007/s11071-005-1914-0